# Semi-parametric second-order efficient estimation of the period of a signal

I. CASTILLO

[1]*Department of Mathematics, Vrije Universiteit Amsterdam, De Boelelaan 1081a, 1081 HV Amsterdam, The Netherlands. E-mail: i.castillo@few.vu.nl*

This paper is concerned with the estimation of the period of an unknown periodic function in Gaussian white noise. A class of estimators of the period is constructed by means of a penalized maximum likelihood method. A second-order asymptotic expansion of the risk of these estimators is obtained. Moreover, the minimax problem for the second-order term is studied and an estimator of the preceding class is shown to be second order efficient.

*Keywords:* exact minimax asymptotics; penalized maximum likelihood; second-order efficiency; semi-parametric estimation; unknown period

## 1. Introduction

*The framework.* Let $x$ be a random process defined by the equation

$$\mathrm{d}x(t) = f(t/\theta)\,\mathrm{d}t + \mathrm{d}W(t), \qquad t \in [-T/2, T/2], \theta > 0,\ T > 0, \tag{1}$$

where $f$ is an unknown real periodic function with period 1, $\theta$ is a period parameter which we seek to estimate and $W(t)$ is the standard Brownian motion on $[-T/2, T/2]$. It is assumed that $x(t)$ is observed continuously on $[-T/2, T/2]$ and that $\theta$ belongs to an interval $\Theta = [\alpha_T, \beta_T]$, where $\alpha_T, \beta_T$ are reals such that $0 < \alpha_T < \beta_T < +\infty$. Throughout the paper, the asymptotic framework $T \to +\infty$ is considered.

*Motivation.* When $f$ is known, under smoothness and identifiability conditions, the maximum likelihood estimator is asymptotically efficient; see Ibragimov and Has'minskii [12], Theorem 5.1, page 203. Here, we are interested in the semi-parametric problem: the parameter $\theta$ must be estimated, but $f$ is an infinite-dimensional nuisance parameter. In a seminal paper, Golubev [9] gives an asymptotically efficient estimator of the period in this framework. In a recent article, Gassiat and Lévy-Leduc [7] obtain an asymptotically efficient estimator of the period in a discretized version of (1) and address the problem of estimation of multiple periods when the signal is a sum of different periodic functions.

From a practical point of view, the problem of period estimation arises in many different areas, such as communications, seismic signal processing and laser vibrometry (see Prenat [16]). Model (1) can be seen as the signal observed by a receptor situated at some distance from a vibrating source. The frequency of the source vibration equals $\theta^{-1}$ and the noise $dW(t)$ stands for the degradation of the original signal due to the distance between the source and receptor. For other practical applications and further references on the subject, we refer to Gassiat and Lévy-Leduc [7] and Lavielle and Lévy-Leduc [14]. In the latter, the authors propose a practical method for estimating the frequency in the discrete-time model cited above, when the number of observations at hand is fixed.

*Semi-parametric estimators and second-order efficiency.* We now return to the general semi-parametric framework. The aim is to estimate the parameter without knowing the nuisance function $f$ and, if possible, to obtain the same optimal asymptotic variance for the estimator as in the parametric framework (see van der Vaart [18], Chapter 25, or Bickel, Klaassen and Wellner [1]). If this is possible, which is the case in model (1), we say that there is no information loss.

However, for a given model, there is often a large choice of asymptotically efficient estimators. This motivates the study of the second order term for the estimator's quadratic risk. This problem was studied for partial linear models by Golubev and Härdle [10]. In this paper, the authors construct second-order efficient estimators when the nuisance function belongs to a known functional class. Golubev and Härdle [11] give nonparametric adaptive versions of their estimators. For essentially nonlinear models, the first result on semi-parametric second-order efficiency was established by Dalalyan, Golubev and Tsybakov [6]. In this paper, the translation model $x(t) = f(t-\theta) + \varepsilon n(t)$, $t \in [-1/2, 1/2]$, is considered and the authors study, as $\varepsilon$ tends to zero, the asymptotic second-order properties of a class of penalized maximum likelihood estimators. Second-order adaptive versions of these results are obtained by Dalalyan [5].

As noted by Dalalyan *et al.* [6], the study of this problem is of particular interest since, *asymptotically*, the second-order terms in a semi-parametric setting are often not much smaller than the first-order terms. Typically, as we shall see for model (1), when the first-order term is $T^{-3/2}$, the second-order term can be of the order of $T^{-3/2} \times T^{-2/5}$.

*Objective and results.* The goal of this paper is to investigate the second-order efficiency in the semi-parametric problem of period estimation. One of the key tools in the paper by Dalalyan *et al.* [6] is the fact that the model can be projected onto a basis and thus described as a discrete collection of independent submodels. As we shall see in the sequel, this property does not hold in our framework. Thus we must work globally on the whole model and hence introduce appropriate methodologies, for instance, the formulation of the criterion upon which the estimation is based as the action of a compact operator on some functional space.

In this paper, we obtain the second-order risk term over a large family of period estimators. We also obtain a lower bound for the second-order term, which is achieved by an estimator of the preceding family. This is the first result about second-order estimation in the period model. It provides a theoretical basis for the choice of the smoothing parameters in the construction of the estimator. It also highlights the important role played by the nonparametric nuisance part of the model. For instance, it originally motivated the work of Castillo, Lévy-Leduc and Matias [3], where the authors study the problem of sharp adaptive estimation of $f$ in model (1).

*Structure of the paper.* In Section 2, we construct a class of estimators using the penalized maximum likelihood method. Section 3 contains our main asymptotic results: we give the second-order properties of the preceding estimators and study the lower bound for the risk over all estimators. Moreover, we construct an estimator achieving the optimal second-order rate. Section 4 is devoted to technical proofs.

## 2. Penalized maximum likelihood estimator

Let us first introduce some useful notation. As the function $f$ in (1) is 1-periodic, we may assume that $f$ can be written as a convergent Fourier series,

$$f(x) = \sum_{k\in\mathbb{Z}} c_k e^{2i\pi kx} = a_1 + \sqrt{2}\sum_{k\geq 1} a_{2k}\cos(2\pi kx) + a_{2k+1}\sin(2\pi kx) = \sum_{k\geq 1} a_k \varepsilon_k(x), \quad (2)$$

where $\varepsilon_1(x) = 1$ and, for $k \geq 1$, $\varepsilon_{2k}(x) = \sqrt{2}\cos(2\pi kx)$ and $\varepsilon_{2k+1}(x) = \sqrt{2}\sin(2\pi kx)$.

The Fisher information for a continuously differentiable $f$ in model (1) is, as $T$ tends to infinity,

$$I_T(f,\theta) = \{1 + o(1)\}\frac{T^3}{12\theta^4}\sum_{k\in\mathbb{Z}}(2\pi k)^2|c_k|^2. \quad (3)$$

Note that this quantity is asymptotically one quarter of the one given by Ibragimov and Has'minskii [12], page 209, since here, the observation interval is $[-T/2, T/2]$ and not $[0,T]$. Also, note that it depends on *both* $f$ and $\theta$. To simplify the notation, we denote it by $I_T$ in the sequel when we do not want to emphasize this dependence.

We now construct a family of estimators of the period which allows us to deal with efficiency at the second order. We use the method of penalized maximum likelihood introduced by Dalalyan *et al.* [6]. Note that here, model (1) cannot be partitioned into a discrete collection of projected one-dimensional submodels since the space of projection would be $L^2([0,\theta])$ and $\theta$ is unknown. Thus the object of our study is the global likelihood function of the model.

The likelihood for estimating $\theta$ in model (1) depends on $\theta$ and $f$, which is given by its Fourier coefficients $(a_k)$. To eliminate the nuisance parameters $(a_k)$, we first assume that the $a_k$'s are independent centered Gaussian random variables with variance $\sigma_k^2$ and independent of the noise, that is, we put a prior distribution on $(a_k)$. To estimate the $a_k$'s, we maximize the posterior distribution of $(a_k)$ given to the observations $\{x(t)\}$. Note that this is equivalent to the maximization of the joint likelihood of $(\{x(t)\},(a_k))$ or of its logarithm, which, thanks to the Girsanov formula (see, e.g., Ibragimov and Has'minskii [12], Appendix 2), is given by the following function $\Phi$:

$$\Phi[\tau,\{x(t)\},(a_k)] = -\frac{1}{2}\int_{-T/2}^{T/2} f(t/\tau)^2\,\mathrm{d}t + \int_{-T/2}^{T/2} f(t/\tau)\,\mathrm{d}x(t) - \sum_{k\geq 1}\frac{a_k^2}{2\sigma_k^2}. \quad (4)$$

Taking the partial derivatives of the last quantity with respect to the Fourier coefficients $a_k$ and using the approximation $\int_{-T/2}^{T/2} \varepsilon_k(t/\tau) f(t/\tau) \, dt \approx a_k T$, we find that the maximum is approximately obtained for $(a_k^*)$ such that

$$(T + \sigma_k^{-2}) a_k^* = \int_{-T/2}^{T/2} \varepsilon_k(t/\tau) \, dx(t).$$

Let $f^*$ be the function with Fourier coefficients $(a_k^*)$. Then $\int_{-T/2}^{T/2} f^*(t/\tau) \, dx(t) = \sum_{k \geq 1} (T + \sigma_k^{-2}) a_k^{*2}$ and $\int_{-T/2}^{T/2} f^*(t/\tau)^2 \, dt \approx T \sum_{k \geq 1} a_k^{*2}$. Thus

$$\Phi[\tau, \{x(t)\}, (a_k^*)] \approx \sum_{k \geq 1} \frac{1}{T + \sigma_k^{-2}} \left( \int_{-T/2}^{T/2} \varepsilon_k(t/\tau) \, dx(t) \right)^2.$$

For reasons of symmetry (for the minimax problem, we will assume that $f$ lies in a Sobolev ellipsoid, see Section 3.3), for $k \geq 1$, we impose the identity $\sigma_{2k} = \sigma_{2k+1}$. This means that we put the same weights on sine and cosine for a given frequency. Hence, writing, for $k \geq 0$, $\lambda_k = (T + \sigma_{2k+1}^{-2})^{-1} T$, which is in $[0, 1]$, we obtain the weighted criterion

$$\mathrm{L}(\tau) = \sum_{k \geq 1} \frac{\lambda_k}{T} \left| \int_{-T/2}^{T/2} e^{2ik\pi t/\tau} \, dx(t) \right|^2. \tag{5}$$

Note that we have dropped the first term in the sum since it does not depend on $\theta$. Moreover, there is no restriction associated with using the Fourier basis; the preceding construction can be used with any orthonormal basis $\{\varepsilon_k\}$ of $L^2[0,1]$.

Equation (5) can be seen as a weighted version of the estimator proposed by Golubev [9]. Note that the weight $\lambda_k$ is outside the square of the integral. We are thus introducing a weight on the 'energy' and not directly on the data, contrary to data tapers methods, which have been studied in the context of frequency estimation (among others) by Chen, Wu and Dahlhaus [4]. In fact, the results obtained by the two methods are of different natures. We also note from equation (4) that the estimation of the nuisance parameters $a_k$ is done by an (approximate) *penalized* maximum likelihood method and thus the 'Bayesian approach' leads, in fact, to penalization (see, e.g., Kimeldorf and Wahba [13]).

We are now able to construct our estimator. As noted by Golubev [9] or Gassiat and Lévy-Leduc [7], direct maximization of (5) would not allow us to distinguish the multiples of the unknown period. To avoid this problem, we aim to take the smallest approximate minimizer of $\mathrm{L}(\tau)$ by choosing

$$\mathcal{E}_T = \left\{ \tau \in \Theta, \ \mathrm{L}(\tau) \geq (1 - \log^{-1/4} T) \sup_{\tau \in \Theta} \mathrm{L}(\tau) \right\}, \tag{6}$$

$$e_T = \inf \mathcal{E}_T. \tag{7}$$

Now, let $B(x,R) = \{\tau \in \Theta, \ |\frac{x}{\tau} - 1| < R\}$. We define our estimator as

$$\theta^* = \underset{\tau \in B(e_T, 1/4)}{\text{Arg max}} \ \mathrm{L}(\tau). \tag{8}$$

In order to understand the behavior of the criterion $\mathrm{L}(\tau)$, we introduce some useful notation. We define symmetrized weights over $\mathbb{Z}$ by letting $\lambda_{-k} = \lambda_k$ and we set $\lambda_0 = 0$. The symbols $\sum_k$, $\sum_{k \neq 0}$, $\sum_{k \geq 0}$ denote sums over $\mathbb{Z}$, $\mathbb{Z}^*$ and $\mathbb{N}$, respectively. Then

$$\mathrm{L}(\tau) = [\Gamma(\tau) + \mathrm{X}(\tau) + \Psi(\tau)]/2,$$

where

$$\begin{cases} \Gamma(\tau) = \sum_k \lambda_k T^{-1} \left| \int_{-T/2}^{T/2} e^{2ik\pi t/\tau} f(t/\theta) \, dt \right|^2, \\ \mathrm{X}(\tau) = 2 \sum_k \lambda_k T^{-1} \int_{-T/2}^{T/2} e^{2ik\pi t/\tau} f(t/\theta) \, dt \int_{-T/2}^{T/2} e^{-2ik\pi t/\tau} \, dW(t), \\ \Psi(\tau) = \sum_k \lambda_k T^{-1} \left| \int_{-T/2}^{T/2} e^{2ik\pi t/\tau} \, dW(t) \right|^2. \end{cases} \tag{9}$$

## 3. Second-order asymptotics in the period model

Let us begin this section with some definitions. We say that a function is $o(1)$ (resp., $O(1)$) if it tends to zero (resp., is bounded) as $T$ goes to infinity. If $(z_k)$ is a sequence of complex numbers indexed by $\mathbb{Z}$, then

$$\|z\|^2 = \sum_k |z_k|^2, \qquad \|z\|_1 = \sum_k |z_k| \quad \text{and, for } m \geq 1, \quad \|z^{(m)}\|^2 = \sum_k |z_k|^2 (2\pi k)^{2m}.$$

We write simply $\|z'\|$ instead of $\|z^{(1)}\|$. If $v$ is a square-integrable function, $\|v\|$ denotes its $\mathrm{L}^2$-norm on $[0,1]$. In the sequel, $C$ denotes a universal constant.

### 3.1. Assumptions on the model

Let us recall that $\Theta = [\alpha_T, \beta_T]$. We assume that as $T$ tends to $+\infty$,

$$(\mathbf{P1}) \quad \alpha_T^{-1} = O(T) \qquad (\mathbf{P2}) \quad \beta_T = O(\log T).$$

Let us assume that $f$ belongs to some class $F = F(\rho, C_0)$ of smooth functions whose Fourier coefficients $(c_k)$ satisfy the Fourier expansion (2) and

$$(\mathbf{F1}) \quad |c_1|^2 \geq \rho > 0 \qquad (\mathbf{F2}) \quad \sum_k (2\pi k)^4 |c_k|^2 \leq C_0 < +\infty.$$

We consider sequences of weights $(\lambda_k)$ such that $\lambda_0 = 0$, $\lambda_1 = 1$ and for all integers $k$, $\lambda_{-k} = \lambda_k$ and $0 \leq \lambda_k \leq 1$. We also assume that there exists $N_T$ tending to $+\infty$ such that

$$(\mathbf{W0}) \quad \lambda_k = 0 \quad \text{for } k \geq N_T \quad \text{and} \quad N_T^4 = o(T) \quad \text{as } T \text{ tends to } +\infty.$$

Moreover, we assume that there exist positive constants $\rho_1$ and $C_1$ such that

$$(\mathbf{W1}) \quad \|\lambda'\| \geq \rho_1 \log^2 T \max_{k \geq 1} \lambda_k (2\pi k);$$

$$(\mathbf{W2}) \quad \sum_k \lambda_k (2\pi k)^4 \leq C_1 T.$$

Finally, we use the following technical assumption:

$$(\mathbf{T}) \quad \left[\sum_k (1-\lambda_k)(2\pi k)^2 |c_k|^2\right]^2 = o\left[\frac{1}{\log T}\sum_k (1-\lambda_k)^2 (2\pi k)^2 |c_k|^2\right].$$

The $o$ and $O$ in the previous assumptions are meant to be *uniform* with respect to $f$ and $\lambda$. In the sequel, all $o$ and $O$ will be. Assumptions (**F1**) and (**F2**) are regularity conditions on $f$. Assumptions (**W0**), (**W1**) and (**W2**) are satisfied for quite a large variety of weight sequences. For instance, they are fulfilled for projection weights $(\mathbf{1}_{|k| \leq N_T})$, provided that $N_T \geq C \log^4 T$ and $N_T^5 = O(T)$.

Let us briefly study two consequences of the preceding assumptions. First, (**F1**) and (**F2**) imply that there exists a constant $h > 0$—for example $h = \rho/C_0$—such that for any integer $p \geq 2$,

$$\sum_{q \neq 0} |c_{pq}|^2 \leq (1-h) \sum_{q \neq 0} |c_q|^2. \tag{10}$$

Second, since the weights are between 0 and 1, we have $(1-\lambda_k)^2 \leq (1-\lambda_k)$ and thus

$$\sum_k (1-\lambda_k)(2\pi k)^2 |c_k|^2 \times \sum_k (1-\lambda_k)^2 (2\pi k)^2 |c_k|^2 \leq \left[\sum_k (1-\lambda_k)(2\pi k)^2 |c_k|^2\right]^2.$$

Using (**T**), we obtain

$$\sum_k (1-\lambda_k)(2\pi k)^2 |c_k|^2 = o(\log^{-1} T). \tag{11}$$

In particular, note that for each fixed $k$, the weight $\lambda_k$ tends to 1 as $T$ tends to $+\infty$.

### 3.2. Second-order asymptotics for the risk

Let us introduce the following functional, where the $c_k$'s are the Fourier coefficients of $f$.

$$R_T(f, \lambda) = \sum_k (2\pi k)^2 \left((1-\lambda_k)^2 |c_k|^2 + \frac{1}{T}\lambda_k^2\right). \tag{12}$$

This functional corresponds to a term of *nonparametric* estimation. In fact, suppose that we have at hand a Gaussian sequence model defined by $y_k = c_k + T^{-1/2}\zeta_k$, where $\zeta_k = (1/\sqrt{2})(\zeta_{1,k} + i\zeta_{2,k})$ and $\{\zeta_{1,k}\}, \{\zeta_{2,k}\}$ are independent sequences of standard normal random variables. Denoting by $f$ the function with Fourier coefficients $c_k$, suppose we seek to estimate the derivative $f'$ given by the Fourier coefficients $c'_k = 2ik\pi c_k$. Now, consider the linear estimator $\widehat{f'}$ defined by the Fourier coefficients $(\widehat{c'_k}) = (\lambda_k(2ik\pi)y_k)$. Then (12) is nothing but the quadratic risk $\mathbf{E}(\|\widehat{f'} - f'\|^2)$.

**Theorem 1.** *Under assumptions* (**P**), (**F**), (**W**) *and* (**T**), *uniformly in* $f \in F$ *and in* $\theta \in \Theta$, *the estimator* $\theta^*$ *defined by* (8) *satisfies, as* $T$ *tends to* $+\infty$,

$$\mathbf{E}_{\theta,f}((\theta^* - \theta)^2 \mathrm{I}_T(f,\theta)) = 1 + \{1 + o(1)\}\frac{R_T(f,\lambda)}{\|f'\|^2}. \tag{13}$$

*Remark 1.* The preceding theorem gives a second-order expansion for the risk. By (11), $\sum_k (2\pi k)^2 (1 - \lambda_k)^2 |c_k|^2 = o(\log^{-1} T)$. Moreover, (**W0**) gives $T^{-1}\sum_k (2\pi k)^2 \lambda_k^2 \leq CT^{-1}N_T^3 = o(T^{-1/4})$. Hence, $R_T(f,\lambda) = o(1)$ as $T$ tends to $+\infty$.

*Remark 2.* The second-order term $\|f'\|^{-2}R_T(f,\lambda)$ is the same as the one obtained by Dalalyan *et al.* [6] and similar to the one obtained by Golubev and Härdle [10]. It seems to be a general feature of smooth semi-parametric models that the term of order 2 reflects the estimation of the score function of the model.

*Remark 3.* In this paper, we consider model (1) on the centered interval $[-T/2, T/2]$. It can be checked that for a different interval of length $T$, for example $[0,T]$, equation (13) still holds.

**Proof of Theorem 1.** Here, we give the main lines of the proof, the technical aspects being dealt with in Section 4. A transversal object is the criterion L on which the definition of $\theta^*$ relies. We study its deterministic part $\Gamma$ in Section 4.1 by means of Fourier techniques. In Section 4.2, we study the stochastic parts X and $\Psi$ by introducing a well-chosen operator on $\mathrm{L}^2[-T/2, T/2]$.

The first step is to obtain consistency results for $\theta^*$. Let us define the event

$$\mathcal{A}_0 = \left\{\sup_{\tau \in \Theta} |\mathrm{X}(\tau) + \Psi(\tau)| \leq T^{3/4}\right\}. \tag{14}$$

The complement of $\mathcal{A}_0$ has probability $o(T^{-p})$ for any integer $p$, by Lemma 6.

Due to Lemma 1, with $\gamma_T = T^{1/8}$, we have $\sup_{\tau \in \Theta} \Gamma(\tau) \leq T\{\sum_{k \neq 0} |c_k|^2 + o(\gamma_T^{-1/2})\}$. Thus, on the event $\mathcal{A}_0$, we have $\sup_{\tau \in \Theta} \mathrm{L}(\tau) \leq T\{\sum_{k \neq 0} |c_k|^2 + o(\gamma_T^{-1/2}) + T^{-1/4}\}/2$.

On the other hand, again using Lemma 1 and the definition of $\mathcal{A}_0$, we have, on $\mathcal{A}_0$, $2L(\theta) \geq T\{\sum_{k \neq 0} |c_k|^2 + o(\log^{-1} T) - T^{-1/4}\}$. Thus $\theta \in \mathcal{E}_T$ on $\mathcal{A}_0$ for sufficiently large $T$. For any integers $p \geq 2$ and $j \geq 1$, bounding the weights from above by 1 and using (10), we obtain $\sum_q \lambda_{qj}|c_{qp}|^2 T \leq \sum_{q \neq 0}|c_{qp}|^2 T \leq (1-h)\sum_{q \neq 0}|c_q|^2 T$. Hence, with $\varepsilon_T = T^{-3/8}$, we have $\mathcal{E}_T \subset \bigcup_{j \geq 1} B(j\theta, \varepsilon_T)$ on $\mathcal{A}_0$. Therefore, $e_T \in B(\theta, \varepsilon_T)$ on $\mathcal{A}_0$.

The definition of $\theta^*$ then implies that $\theta^* \notin B(2\theta, \varepsilon_T)$ on $\mathcal{A}_0$, thus

$$\theta^* \in B(\theta, \varepsilon_T) \qquad \text{on } \mathcal{A}_0.$$

In fact, we need to refine this rate. For any $D > 0$, let us define the event $\mathcal{A}_1$ as

$$\mathcal{A}_1 = \{|\theta^* - \theta|T^{3/2}\theta^{-2} \leq D\log^{1/2} T\}. \tag{15}$$

Lemma 11 establishes that for every integer $p$, for sufficiently large $D$, $\mathbf{P}(\mathcal{A}_1^c)$ is $o(T^{-p})$.

The next step is to introduce the variable $\widehat{\tau}$, which is a theoretical tool for our proof (of course, it is not an estimator since it requires knowledge of $\theta$):

$$\mathrm{L}'(\theta) + (\widehat{\tau} - \theta)\mathbf{E}(\mathrm{L}''(\theta)) = 0. \tag{16}$$

By Lemma 12, $\widehat{\tau}$ has the desired second-order expansion. It remains to show that $\theta^*$ and $\widehat{\tau}$ are sufficiently close, that $\theta^*$ has the same expansion as $\widehat{\tau}$ at order 2.

$$\mathbf{E}((\theta^* - \theta)^2 \mathrm{I}_T) = \mathbf{E}((\theta^* - \theta)^2 \mathrm{I}_T \mathbf{1}_{\mathcal{A}_1}) + \mathbf{E}((\theta^* - \theta)^2 \mathrm{I}_T \mathbf{1}_{\mathcal{A}_1^c}).$$

Using (**P1**), (**P2**) and (43), we arrive at the conclusion that

$$\mathbf{E}((\theta^* - \theta)^2 \mathrm{I}_T \mathbf{1}_{\mathcal{A}_1^c}) \leq C\beta_T^2 \alpha_T^{-4} T^3 \mathbf{P}(\mathcal{A}_1^c),$$

which is negligible compared to a given power of $T$ if $D$ is chosen sufficiently large. Now,

$$\mathbf{E}((\theta^* - \theta)^2 \mathrm{I}_T \mathbf{1}_{\mathcal{A}_1})$$
$$= \mathbf{E}((\widehat{\tau} - \theta)^2 \mathrm{I}_T) + \mathbf{E}([(\theta^* - \widehat{\tau})^2 + 2(\theta^* - \widehat{\tau})(\widehat{\tau} - \theta)]\mathrm{I}_T \mathbf{1}_{\mathcal{A}_1}) - \mathbf{E}((\widehat{\tau} - \theta)^2 \mathrm{I}_T \mathbf{1}_{\mathcal{A}_1^c}).$$

The first term in the above expression gives the appropriate expansion, by Lemma 12, while the last one is negligible due to (43), (**P1**) and (**P2**). The middle term must still be dealt with and this is carried out by Lemma 13. Hence, the expansions of the risk for $\theta^*$ and $\widehat{\tau}$ are the same, which proves Theorem 1. $\square$

### 3.3. Lower bound for the risk

In this section, we establish a lower bound, in the minimax sense, for the second-order term over all possible estimators of $\theta$. For any $\beta > 0$, $L > 0$, let us define

$$\mathcal{W}_{\beta,L} = \left\{ f = \{c_k\}_{k \in \mathbb{Z}}, \sum_k |2\pi k|^{2\beta}|c_k|^2 \leq L \right\}.$$

The expression of the second-order term in Theorem 1 suggests minimizing the functional $R_T(f, \lambda)$. The behavior of this functional is well understood; see Pinsker [15] or Tsybakov [17], Chapter 3 for a complete overview of the subject. There exists a function $s$ in $\mathcal{W}_{\beta,L}$ and a sequence $q$ in $l^2$ such that $(s,q)$ is a saddle point of $R_T(f,\lambda)$ over $\mathcal{W}_{\beta,L} \times l^2(\mathbb{Z})$. We have

$$r_T = R_T(s,q) = \inf_{\lambda \in l^2} \sup_{f \in \mathcal{W}_{\beta,L}} R_T(f,\lambda) = \frac{1}{T}\sum_k (2\pi k)^2 q_k. \tag{17}$$

We also have explicit expressions for $s$ and $q$ in terms of the solution $W_T$ of the equation

$$\frac{1}{T}\sum_{k \in \mathbb{Z}^*} \left[\left(\frac{W_T}{|k|}\right)^{\beta-1} - 1\right]_+ (2\pi|k|)^{2\beta} = L. \tag{18}$$

Let $g$ be a function in $F(\rho, C_0)$. Let us denote by $g_k$ its Fourier coefficients. For any $\delta > 0$ and fixed $\beta > 0, L > 0$, we define a neighborhood of $g$ as follows:

$$F_\delta(g) = \{f = g + v, \|v\| \leq \delta, v \in \mathcal{W}_{\beta,L}\}.$$

**Theorem 2.** *Suppose that $\beta \geq 2$ and $L > 0$. For any $\delta_T \to 0$ such that there exists $\alpha > 0$ with $\delta_T^2 T/W_T^{1+\alpha} \to +\infty$, as $T$ tends to $+\infty$,*

$$\inf_{\hat\theta} \sup_{\theta, f \in F_{\delta_T}(g)} \mathbf{E}_{\theta,f}((\hat\theta - \theta)^2 I_T(\theta,f)) \geq 1 + \{1 + o(1)\}\frac{r_T}{\|g'\|^2},$$

*where the infimum is taken over all estimators $\hat\theta$ based on the observations $\{x(t)\}$ and where $r_T$ is defined by* (17).

**Proof.** The first step is to change variables so that, contrary to (3), the Fisher information will no longer depend on the parameter of interest. Let us define $\omega = \theta^{-1}$. It is not hard to check that the Fisher information $J_T(f)$ for estimating $\omega$ in model (1) equals

$$J_T(f) = \theta^4 I_T(\theta,f) = \{1 + o(1)\}T^3\|f'\|^2/12$$

and, in particular, up to the $o(1)$ term, does not depend on $\theta$. Since the mapping $\hat\theta \to \hat\theta/\theta^2$ between the set of estimators of $\theta$ and the set of estimators of $\omega$ is one to one, we obtain

$$\inf_{\hat\theta} \sup_{\theta \in \Theta, f \in F_{\delta_T}(g)} \mathbf{E}_{\theta,f}((\hat\theta - \theta)^2 I_T(\theta,f)) = \inf_{\hat\omega} \sup_{\omega, f \in F_{\delta_T}(g)} \mathbf{E}_{\omega,f}((\hat\omega - \omega)^2 J_T(f)),$$

where the last infimum is taken over all estimators $\hat\omega$ based on the observations $\{x(t)\}$ and where the supremum is taken over $\omega \in [\beta_T^{-1}, \alpha_T^{-1}]$.

The second step is to bound the minimax risk from below by a well-chosen Bayes risk. Let us define a prior on $f$ putting a prior on each Fourier coefficient of $f$, concentrating around the Fourier coefficient $g_k$ of $g$ in the following way:

$$f_k = g_k + \sigma_k\{\xi_{k,1} + i\xi_{k,2}\},$$

where for $p = 1, 2$, $\xi_{k,p} \sim \mathcal{N}(0, 1/2)$ are independent and independent of $\{W(t)\}$ and

$$\sigma_k^2 = \begin{cases} 0, & |k| \leq \gamma_T W_T, \\ (1 - \gamma_T)|s_k|^2, & |k| > \gamma_T W_T, \end{cases} \quad (19)$$

where $\gamma_T = 1/\log T$ and the $s_k$'s are the Fourier coefficients of the function defined in (17). We then choose as prior for $\omega$ a distribution with density $\pi$ on $[\alpha_T, \beta_T]$, where $\pi$ vanishes at the endpoints of the interval and such that the Fisher information $J_\pi = \int \pi'(\omega)^2 \pi^{-1}(\omega) \, dx$ is finite.

We denote by $\Psi_\sigma(f)$ the distribution associated with the preceding prior on $f$. We then let

$$\overline{J_T} = \int J_T(f) \, d\Psi_\sigma(f) = \{1 + o(1)\} \frac{T^3}{12} \sum_k (2\pi k)^2 (|g_k|^2 + 2\sigma_k^2).$$

In the sequel, we shall denote by $\mathbf{E}$ the expectation in the full Bayesian model.

For the preceding choice of $\sigma_k^2$, the random function $f$ is close to $g$ with high probability. More precisely, one can check (as in Dalalyan *et al.* [6]) that $\mathbf{P}(F_{\delta_T}(g)^c)$ decreases at an exponential rate. Moreover, it is also not difficult to check that with this choice of $\sigma_k^2$, we have

$$|12\overline{J_T} T^{-3} - \|g'\|^2| = o(1). \quad (20)$$

$$T^{-1} \sum_k (2\pi k)^2 \mu_k = \{1 + o(1)\} r_T, \qquad \text{with } \mu_k = \frac{\sigma_k^2}{T^{-1} + \sigma_k^2}. \quad (21)$$

Let us bound the minimax risk below using the full Bayesian model:

$$\inf_{\hat{\omega}} \sup_{\omega, f \in F_{\delta_T}(g)} \mathbf{E}_{\omega, f}((\hat{\omega} - \omega)^2 J_T(f)) \geq \inf_{\hat{\omega}} \mathbf{E}((\hat{\omega} - \omega)^2 \overline{J_T} \mathbf{1}_{F_{\delta_T}(g)}(f)) \quad \text{(I)}$$

$$- \sup_{\hat{\omega}} \mathbf{E}((\hat{\omega} - \omega)^2 (\overline{J_T} - J_T(f)) \mathbf{1}_{F_{\delta_T}(g)}). \quad \text{(II)}$$

The outline of the proof of the theorem is then as follows. We bound (I) from below:

$$\text{(I)} \geq \inf_{\hat{\omega}} \overline{J_T} \mathbf{E}((\hat{\omega} - \omega)^2) - C\beta_T^2 T^3 \mathbf{P}(F_{\delta_T}(g)^c).$$

The second term in this difference is negligible, due to the exponential bound on $\mathbf{P}(F_{\delta_T}(g)^c)$. Let us denote by $\mathcal{J}_T(\omega)$ the Fisher information associated with the observations $\{x(t)\}$ in the Bayesian model with respect to $f$ *with fixed* $\omega$. The Van Trees inequality (see Gill and Levit [8]) applied to the full Bayesian model gives

$$\inf_{\hat{\omega}} \mathbf{E}((\hat{\omega} - \omega)^2) \geq \frac{1}{\int \mathcal{J}_T(\omega) \pi(d\omega) + J_\pi}.$$

The key point of the proof is then to obtain an expansion of $\mathcal{J}_T(\omega)$. Note that as the model cannot be partitioned into one-dimensional submodels, the Fisher information $\mathcal{J}_T(\omega)$ is not a sum of Fisher information over submodels and must thus has to be handled globally. As $T$ tends to $+\infty$, with the $\sigma_k^2$'s defined by (19), we have

$$\mathcal{J}_T(\omega) = \overline{J_T} - \{1 + o(1)\}\frac{T^2}{12}\sum_k (2\pi k)^2 \mu_k. \tag{22}$$

The proof of (22) is not difficult, but it is technical. For detailed calculations, we refer the reader to Castillo [2].

Using equation (22), we may deduce

$$\inf_{\hat{\omega}} \overline{J_T} \mathbf{E}((\hat{\omega} - \omega)^2) \geq \frac{\overline{J_T}}{\overline{J_T} - \{1 + o(1)\}(T^2/12)\sum_k (2\pi k)^2 \mu_k + J_\pi}.$$

Now, using (20) and (21), we obtain

$$\inf_{\hat{\omega}} \overline{J_T} \mathbf{E}((\hat{\omega} - \omega)^2) \geq \frac{1}{1 - \{1 + o(1)\}\|g'\|^{-2}(r_T - 12T^{-3}J_\pi)}.$$

Using the fact that $T^{-3} = o(r_T)$ and the inequality $(1+x)^{-1} \geq 1 - x$, valid for $x > -1$,

$$\inf_{\hat{\omega}} \overline{J_T} \mathbf{E}((\hat{\omega} - \omega)^2) \geq 1 + \{1 + o(1)\}\|g'\|^{-2} r_T.$$

Finally, (II) is negligible with respect to (I), as in Dalalyan *et al.* [6]. □

*Remark 4.* Changing variables works well in our framework since the Fisher information (3) has a quite simple separated form in $f$ and $\theta$. This might not be the case for other models. For instance, in the Cox model, the form of the efficient Fisher information as given in van der Vaart [18], page 416 is more involved and a change of variables does not seem to make the parameter of interest vanish. We think that this difficulty could be overcome by also carefully choosing a prior on the $\theta$-parameter.

### 3.4. Achieving the lower bound

**Theorem 3.** *Assume that the conditions of Theorem 2 are fulfilled and that there is some $p > \beta \geq 2$ such that the sum $\sum |2\pi k|^{2p}|g_k|^2$ is finite. Then there exists a sequence of weights $\lambda^*$ such that the corresponding estimator $\widehat{\theta}(\lambda^*)$ defined by (8) achieves the bound of Theorem 2. That is, as $T$ tends to $+\infty$,*

$$\sup_{\theta \in \Theta, f \in F_{\delta_T}(g)} \mathbf{E}_{\theta,f}((\widehat{\theta}(\lambda^*) - \theta)^2 I_T(\theta, f)) = 1 + \{1 + o(1)\}\|g'\|^{-2} r_T.$$

Let us comment on the rate of convergence for the second-order term obtained above. Pinsker's theory (see Pinsker [15]) states that, up to a constant which can be computed in terms of $\beta$ and $L$, $r_T$ in [17] is of the order $T^{(2-2\beta)/(2\beta+1)}$. For example, in the case $\beta = 2$, the optimal second-order term is of the order $T^{-2/5}$. This semi-parametric rate is fairly slow compared to the first-order rate, which is, up to a constant, $T^{-3/2}$, highlighting the importance of second-order terms in semi-parametric estimation. Moreover, the less regular the function, the more significant the second-order term.

**Proof of Theorem 3.** The idea is to use slightly modified Pinsker weights for $\lambda^*$ by letting

$$\lambda_k^* = \begin{cases} 1, & |k| \leq \gamma_T W_T, \\ [1 - (|k|/W_T)^{\beta-1}]_+, & |k| > \gamma_T W_T, \end{cases} \quad (23)$$

where $W_T$ satisfies (18). It is not difficult to check that these weights satisfy assumptions (**W**) and (**T**). Due to Theorem 1, it suffices to check that the supremum of $R_T(f, \lambda^*)$ for $f$ in the vanishing neighborhood under consideration is indeed equivalent to $r_T$. There is no difficulty with respect to Dalalyan *et al.* [6], so this argument is omitted. □

## 4. Proof of Theorem 1

**Definition 1.** *We say that the probability of a measurable set $\mathcal{A}$ is negligible if, for all integers $p$, we have $\mathbf{P}(\mathcal{A}(K)) = o(T^{-p})$ as $T$ tends to $+\infty$.*

### 4.1. Behavior of the deterministic part $\Gamma$

Let us denote by $\hat{\phi}$ the Fourier transform of the indicator function of $[-1/2, 1/2]$: $\hat{\phi}(x) = \int_{-1/2}^{1/2} e^{2i\pi x t} dt = \sin(\pi x)/(\pi x)$. In the sequel, we use the following bounds on $\hat{\phi}$. For any $p \in \mathbb{N}$, there exist constants $M_p > 0$, depending only on $p$, such that

$$|\hat{\phi}^{(p)}(u)| \leq M_p \quad \text{for } u \in [-1,1], \qquad |\hat{\phi}^{(p)}(u)| \leq M_p/|u| \quad \text{for } |u| > 1/4. \quad (24)$$

$$|\hat{\phi}'(u)| \leq C_1|u| \quad \text{for all } u \in \mathbb{R}, \qquad |\hat{\phi}^{(3)}(u)| \leq C_2|u| \quad \text{for all } u \in \mathbb{R}. \quad (25)$$

For any real number $x$, let us denote by $\Delta(x)$ its distance to $\mathbb{Z}$ and by $]x[$ the (smallest) integer realizing this distance. Let us also introduce the auxiliary notation

$$a_{k,l} = \frac{T}{\theta}\left(l - \frac{k\theta}{\tau}\right), \qquad b_{k,l} = \frac{T}{\theta}(l-k).$$

Using (2), we obtain

$$\Gamma(\tau) = \sum_k \lambda_k T \left| \sum_l c_l \hat{\phi}(a_{k,l}) \right|^2, \quad (26)$$

$$\Gamma(\tau) = \sum_k \lambda_k T \left[ \left| c_{]k\theta/\tau[} \hat{\phi}\left(\frac{T}{\theta}\Delta(k\theta/\tau)\right) + \sum_{l \neq ]k\theta/\tau[} c_l \hat{\phi}(a_{k,l}) \right|^2 \right]. \quad (27)$$

For any $p, j \in \mathbb{Z}$, we denote by $p \wedge j$ the greatest common divisor of $p$ and $j$.

**Lemma 1.** *Let $\gamma_T = T^{1/8}$ and $\varepsilon_T = T^{-3/8}$. Then, as $T$ tends to $+\infty$,*

$$\Gamma(\tau) = o(\gamma_T^{-1/2} T), \qquad \text{if } \tau \notin \bigcup_{j \geq 1, 0 < p \leq \gamma_T} B(j\theta/p, \varepsilon_T), \tag{28}$$

$$\Gamma(\tau) \leq \left\{ \sum_q \lambda_{qj} |c_{qp}|^2 + o(\gamma_T^{-1/2}) \right\} T, \quad \text{if } \tau \in B(j\theta/p, \varepsilon_T) \ (p \leq \gamma_T, p \wedge j = 1), \tag{29}$$

$$\Gamma(\theta) = \left\{ \sum_{k \neq 0} |c_k|^2 + o(\log^{-1} T) \right\} T. \tag{30}$$

**Proof.** First, note that, without restriction, we can just study the sum over $|l| \leq \gamma_T$ in (26). Indeed, using the Cauchy–Schwarz inequality, $\|\lambda\|_1 \leq N_T = o(T^{1/4})$ and (**F2**), we have

$$\sum_k \lambda_k T \left| \sum_{|l| > \gamma_T} c_l \hat{\phi}(a_{k,l}) \right|^2 \leq \|\lambda\|_1 T \sum_{|l| > \gamma_T} \frac{1}{l^4} \sum_{|l| > \gamma_T} l^4 |c_l|^2 = o(\gamma_T^{-1} T). \tag{31}$$

In the remainder of the proof, $p$ and $j$ are two relatively prime integers such that $p \leq \gamma_T$. Also, note that by symmetry, we can always assume that $k$ is positive.

Assume that $\tau$ is not in a ball $B(j\theta/p, \varepsilon_T)$ with $p \leq \gamma_T$ and $j \geq 1$. For any $l$ such that $|l| \leq \gamma_T$ and any integer $k$, this implies that $|a_{k,l}| \geq T\beta_T^{-1}\varepsilon_T \geq 1$. By (24), we have $|\hat{\phi}(a_{k,l})|^2 \leq \beta_T^2 T^{-2} \varepsilon_T^{-2}$. Using (26) and (31), this yields $\Gamma(\tau) = o(\gamma_T^{-1/2} T)$.

Assume that $\tau$ is in a ball $B(j\theta/p, \varepsilon_T)$ with $p \leq \gamma_T$ and $p \wedge j = 1$. In this case, we observe that the sum over $l \neq ]k\theta/\tau[$ in (27) is negligible, due to the triangle inequality, (**F2**) and (24). There are then two cases for the integer $k$ in (27):

- There exists an integer $q$ such that $k = qj$. Then, for sufficiently large $T$, $]k\theta/\tau[ = pq$ since

$$|k\theta/\tau - pq| < pq\varepsilon_T \leq \gamma_T N_T \varepsilon_T = o(1).$$

- The integer $k$ is not a multiple of $j$. Then $k = qj + r$ with $r < j$. Note that

$$\left| \frac{k\theta}{\tau} - \left(pq + \frac{rp}{j}\right) \right| \leq pq\varepsilon_T + \frac{rp\varepsilon_T}{j} \leq Cpq\varepsilon_T \leq C \frac{\gamma_T N_T \varepsilon_T}{j}.$$

Since $p \wedge j = 1$, $k\theta/\tau$ is at a distance to the integers larger than $1/2j$. Since $\tau$ must lie in $\Theta$, we have $j\theta/p \leq 2\beta_T$. Thus, $j^{-1} \geq \theta/(2p\beta_T)$. Therefore,

$$\left| \frac{T}{\theta} \Delta(k\theta/\tau) \right| \geq \frac{T}{4\gamma_T \beta_T},$$

which allows us to prove that the corresponding term is negligible, using (24).

Finally, writing (27) at $\tau = \theta$, it is easy to check, by similar arguments, that

$$\Gamma(\theta) = \sum_k \lambda_k |c_k|^2 T + o\{(\beta_T T^{-1} + N_T \beta_T^2 T^{-2})T\}.$$

Moreover, due to (11), $\sum_{k \neq 0}(1-\lambda_k)|c_k|^2 \leq \sum_k (1-\lambda_k)(2\pi k)^2 |c_k|^2 = o(\log^{-1} T)$, which establishes (30). □

Let us define the following weighted Fisher information:

$$\mathrm{I}(\lambda) = \frac{T^3}{12\theta^4} \sum_k \lambda_k (2\pi k)^2 |c_k|^2, \qquad \mathrm{I}(\lambda^2) = \frac{T^3}{12\theta^4} \sum_k \lambda_k^2 (2\pi k)^2 |c_k|^2. \tag{32}$$

**Lemma 2.** *As $T$ tends to $+\infty$,*

$$\Gamma'(\theta)^2 = T^3 \theta^{-4} o(\|\lambda'\|^2 T^{-1}), \qquad \Gamma''(\theta) = -2\mathrm{I}(\lambda) + T^3 \theta^{-4} o(\|\lambda'\|^2 T^{-1}).$$

**Proof.** For any $\tau$ in $\Theta$, taking the derivative with respect to $\tau$ in (9),

$$\Gamma'(\tau) = \frac{T^2}{\tau^2} \sum_k \lambda_k k \sum_{p,l} \overline{c_p} c_l (\hat{\phi}'(a_{k,l})\hat{\phi}(a_{k,p}) + \hat{\phi}'(a_{k,p})\hat{\phi}(a_{k,l})). \tag{33}$$

Now, consider the case $\tau = \theta$ and write the sum over $p, l$ in (33) as

$$\sum_{p,l} = \sum_{p=k,l=k} + \left\{ \sum_{p=k,l\neq k} + \sum_{p\neq k, l=k} \right\} + \sum_{p\neq k, l\neq k} = a + b + c. \tag{34}$$

Then note that $a = 0$, $|b| \leq |c_k| \sum_p |c_p| \theta T^{-1}$ and $|c| \leq (\sum_p |c_p|)^2 \theta^2 T^{-2}$. Hence, $\Gamma'(\theta) = O(T/\theta)$, which, using (**P2**) and (**W1**), gives the first expansion. The expansion for $\Gamma''(\theta)$ is obtained similarly. □

Let us define the set $\mathcal{V}_{\mathcal{A}_1}$ as

$$\mathcal{V}_{\mathcal{A}_1} = \{\tau \in \Theta, \ |\tau - \theta| \leq DT^{-3/2} \theta^2 \log^{1/2} T\}. \tag{35}$$

**Lemma 3.** *As $T$ tends to $+\infty$, uniformly in $\tau \in \mathcal{V}_{\mathcal{A}_1}$, we have*

$$\Gamma'(\tau) = O(\theta^{-2} T^{3/2} \log^{1/2} T), \qquad \Gamma''(\tau) = O(\theta^{-4} T^3), \qquad \Gamma^{(3)}(\tau) = o(\theta^{-6} T^4).$$

**Proof.** Using (25) and the fact that $\tau$ lies in $\mathcal{V}_{\mathcal{A}_1}$, $|\hat{\phi}'(a_{k,k})| \leq C_1 |a_{k,k}| \leq CkT^{-1/2} \log^{1/2} T$. We again write the sum (33) as in (34). By to (25) and (**F2**), the term corresponding to $\sum_{p=k=l}$ is bounded above by $CT^{3/2} \theta^{-2} \log^{1/2} T \sum \lambda_k k^2 |c_k|^2 \leq CT^{3/2} \theta^{-2} \log^{1/2} T$. Similarly, one concludes that the term $\sum_{p=k,l\neq k} + \sum_{p\neq k,l=k}$ is bounded above by $C(T\theta^{-1} + T^{3/2} \theta^{-2} \log^{1/2} T)$ and that the term $\sum_{p\neq k, l\neq k}$ is bounded above by $\sum \lambda_k k \leq$

$N_T^2$, which together with (**P2**) and (**W0**), gives the result. The other two expansions are obtained similarly. □

## 4.2. Behavior of the stochastic parts X and $\Psi$

Let $\mathbb{K}_T^\tau$ be the operator on $L^2([-T/2, T/2])$ such that for any $g$ in $L^2([-T/2, T/2])$,

$$\mathbb{K}_T^\tau g(t) = \int_{-T/2}^{T/2} \left\{ T^{-1} \sum_k \lambda_k e^{2i\pi k(t-u)/\tau} \right\} g(u) \, du = \int_{-T/2}^{T/2} K_T^\tau(t-u) g(u) \, du.$$

Then $X(\tau)$ is nothing but $2 \int \mathbb{K}_T^\tau f(\cdot/\theta)(t) \, dW(t)$. Thus, X and its derivative are centered Gaussian random variables of variance

$$\mathbf{E}(X^{(p)}(\tau)^2) = 4 \int_{-T/2}^{T/2} |[\mathbb{K}_T^\tau f(\cdot/\theta)]^{(p)}(t)|^2 \, dt. \tag{36}$$

On the other hand, note that for every $T$ and $\tau$, the operator $\mathbb{K}_T^\tau$ acting on $L^2[-T/2, T/2]$ is self-adjoint and compact. Namely, it is a convolution operator with bounded kernel $K_T^\tau$. The spectral theorem (see, e.g., Yoshida [19]) then states that it can be characterized by eigenvalues $\{\beta_k\}_{k \in \mathbb{Z}}$ and an orthonormal basis of eigenvectors $\{v_k\}_{k \in \mathbb{Z}}$.

Using expansions over the basis $\{v_k\}$, one can check that there exists a sequence $\{\alpha_k\}_{k \in \mathbb{Z}}$ of independent $\mathcal{N}(0,1)$ random variables such that

$$\Psi(\tau) = \int_{-T/2}^{T/2} \int_{-T/2}^{T/2} K_T^\tau(t-u) \, dW(u) \, dW(t) = \sum_k \beta_k \alpha_k^2. \tag{37}$$

This also holds for the derivatives of $\Psi$, replacing the kernel $K_T^\tau$ by its derivatives. Let us denote by $\beta_k^{(p)}$ the eigenvalues of $\mathbb{K}_T^{\tau(p)}$. It then follows from (37) that for $p \geq 0$,

$$\mathbf{E}[\Psi^{(p)}(\tau)] = \sum_k \beta_k^{(p)} = \int_{-T/2}^{T/2} K_T^{\tau(p)}(0) \, dt, \tag{38}$$

$$\mathbf{Var}[\Psi^{(p)}(\tau)] = 2 \sum_k \beta_k^{(p)2} = 2 \int_{-T/2}^{T/2} \int_{-T/2}^{T/2} |K_T^{\tau(p)}(t-u)|^2 \, dt \, du, \tag{39}$$

where **Var** denotes the variance. These formulas allow the deviations of the process $\Psi^{(p)}$ to be controlled through the study of its Laplace transform.

**Lemma 4.** *There exists $C > 0$ such that, for all $\tau \in \Theta$, $\mathbf{E}(X(\tau)^2) \leq CT$.*

**Proof.** Let $\gamma_k = \sum_l c_l \hat{\phi}(a_{k,l})$. Then

$$\mathbf{E}(X(\tau)^2) = 4T \int_{-1/2}^{1/2} \left| \sum_k \lambda_k \sum_l c_l \hat{\phi}(a_{k,l}) \exp\{2ik\pi tT/\tau\} \right|^2 dt,$$

$$\mathbf{E}(\mathrm{X}(\tau)^2) = 4T \sum_k \lambda_k^2 |\gamma_k|^2 + 4T \sum_{k \neq k'} \lambda_k \lambda_{k'} \gamma_k \overline{\gamma_{k'}} O(\beta_T T^{-1}),$$

where we used the fact that $\int_{-1/2}^{1/2} \exp\{2\mathrm{i}(k-k')\pi t T/\tau\}\,\mathrm{d}t = \mathbf{1}_{k=k'} + \mathbf{1}_{k \neq k'} O(\beta_T T^{-1})$ for all integers $k, k'$. The first term in the last sum is similar to $\Gamma(\tau)$ (see (26), with $\lambda_k^2$ replacing $\lambda_k$) and thus it is a $O(T)$ uniformly in $\tau$ thanks to Lemma 1. The second term is $O(\|\lambda\|_1^2 \beta_T T^{-1})$, which proves the lemma. □

**Lemma 5.** *For all $\tau \in \Theta$, for all positive integers $p$,*

$$\mathbf{E}(\mathrm{X}^{(p)}(\tau)^2) \leq C T^{2p+1} \tau^{-4p} \sum_k \lambda_k^2 (2\pi k)^{2p} (|c_{]k\theta/\tau[}|^2 + O(\beta_T^2 T^{-2})),$$

$$\mathbf{E}(\Psi^{(p)}(\tau)^2) \leq C T^{2p} \tau^{-4p} \|\lambda^{(p)}\|^2.$$

**Proof.** Defining $\xi_k(\tau) = \{\sum_l c_l \hat{\phi}(a_{k,l}(\tau))\}^{(p)}$, the derivative being with respect to $\tau$, we have

$$\mathbf{E}(\mathrm{X}^{(p)}(\tau)^2) = 4T \int_{-1/2}^{1/2} \left| \sum_k \lambda_k e^{2\mathrm{i}\pi kT t/\tau} \xi_k \right|^2 \mathrm{d}t$$

$$\leq 4T \left[ \sum_k \lambda_k^2 |\xi_k|^2 + O(\beta_T T^{-1}) \sum_{k \neq m} \lambda_k \lambda_m \xi_k \overline{\xi_m} \right].$$

Since $|\sum_{k \neq m} \lambda_k \lambda_m \xi_k \overline{\xi_m}| \leq (\sum_k \lambda_k^2 |\xi_k|)^2 \leq N_T \sum_k \lambda_k^2 |\xi_k|^2$, by the Cauchy–Schwarz inequality,

$$\mathbf{E}(\mathrm{X}^{(p)}(\tau)^2) \leq CT \sum_k \lambda_k^2 |\xi_k|^2.$$

One now checks that the predominant term coming from the $p$th derivative while evaluating $\xi_k$ is $T^p \tau^{-2p} k^p \hat{\phi}^{(p)}(a_{k,l})$. Hence,

$$|\xi_k|^2 \leq C T^{2p} \theta^{-4p} k^{2p} \sum_l |c_l \hat{\phi}^{(p)}(a_{k,l})|^2$$

$$\leq C T^{2p} \tau^{-4p} (2\pi k)^{2p} \left[ |c_{]k\theta/\tau[}|^2 + O\left( \left[ \sum_l |c_l| \beta_T T^{-1} \right]^2 \right) \right],$$

which proves the result for $\mathrm{X}^{(p)}$.

To bound $\Psi^{(p)}$, we use (39). Again, one checks that when evaluating the $p$th derivative, there is one dominating term, which is

$$\mathbf{E}(\Psi^{(p)}(\tau)^2) \leq C T^{-2} \tau^{-4p} \sum_{k,l} \lambda_k \lambda_l (2\pi)^{2p} (kl)^p$$

$$\times \int_{-T/2}^{T/2} \int_{-T/2}^{T/2} (t-u)^{2p} \exp\{2\mathrm{i}\pi(t-u)(k-l)/\tau\}\,\mathrm{d}u\,\mathrm{d}t$$

$$\leq CT^{-2}\tau^{-4p} \sum_k \lambda_k^2 (2\pi k)^{2p} O(T^{2p+2}). \qquad \square$$

**Remark 5.** In the sequel, we use Lemma 5 for neighborhoods $\mathcal{V}$ of $\theta$ such that, for any $\tau \in \mathcal{V}$ and any integer $k$ such that $|k| \leq N_T$, we have $]k\theta/\tau[ = k$.

**Lemma 6.** *For any integer $p$, as $T$ tends to $+\infty$,*

$$\mathbf{P}\left(\sup_{\tau \in \Theta} |\mathrm{X}(\tau) + \Psi(\tau)| > \log^2 T\left[T^{1/2} + \sum_k \lambda_k\right]\right) = o(T^{-p}).$$

**Proof.** The proof is not difficult, using Lemma 14 (as in Golubev [9]), and is thus omitted. $\qquad \square$

**Lemma 7.** *If $\mathrm{I}(\lambda^2)$ is defined by (32), then as $T$ tends to $+\infty$,*

$$\mathbf{E}(\mathrm{X}'(\theta)^2) = 4\,\mathrm{I}(\lambda^2) + T^3\theta^{-4}o(\|\lambda'\|^2 T^{-1}).$$

**Proof.** Let $y_s = \sum_l c_l(2\mathrm{i}\pi)\hat{\phi}(b_{l,s})$ and $z_s = -\sum_l c_l\hat{\phi}'(b_{l,s})$. Then

$$T^{-3}\theta^4 \mathbf{E}(\mathrm{X}'(\theta)^2) = 4 \int_{-1/2}^{1/2} \left|\sum_k \lambda_k k(ty_k - z_k) \exp\{2\mathrm{i}\pi ktT/\theta\}\right|^2 \mathrm{d}t$$

$$= 4 \sum_k \lambda_k^2 k^2 (|y_k|^2/12 + |z_k|^2)$$

$$+ O(\theta T^{-1}) \sum_{k \neq p} \lambda_k \lambda_p kp(|y_k \overline{y_p}| + |y_k \overline{z_p}| + |\overline{y_p} z_k| + |z_k \overline{z_p}|).$$

Now, note that $y_s = 2\mathrm{i}\pi c_s + O(\theta T^{-1})$ and $z_s = O(\theta T^{-1})$. This yields the result using, for the remainder term, the fact that, by the Cauchy–Schwarz inequality and (**F2**), we have $\sum_k \lambda_k^2 k^2 |c_k| \leq C\|\lambda\|$ and $\sum_k \lambda_k k |c_k| \leq C$. $\qquad \square$

**Lemma 8.** *As $T$ tends to $+\infty$,*

$$\mathbf{E}(\Psi'(\theta)^2) = T^2\theta^{-4}\{1 + o(1)\} \sum_k (2\pi k)^2 \lambda_k^2/3.$$

**Proof.** Using (36), we obtain

$$\mathbf{E}(\Psi'(\theta)^2) = 8\pi^2 T^2 \theta^{-4} \sum_{k,p} kp\lambda_k \lambda_p \int_{-T/2}^{T/2} \int_{-T/2}^{T/2} \exp(2\mathrm{i}\pi(t-u)(k-p)/\theta)(t-u)^2\,\mathrm{d}t\,\mathrm{d}u,$$

$$\mathbf{E}(\Psi'(\theta)^2) = T^2\theta^{-4}\{(1/3) + O(\theta^2 T^{-2})\} \sum_k (2\pi k)^2 \lambda_k^2.$$

□

**Lemma 9.** *If we recall that $\mathcal{M}_T = \max_k \lambda_k(2\pi k)$, then the quantities*

$$\mathbf{P}\left(\sup_{\mathcal{V}_{\mathcal{A}_1}} |\mathrm{X}''(\tau)| > T^{5/2}\log T\theta^{-4}\right), \qquad \mathbf{P}\left(\sup_{\mathcal{V}_{\mathcal{A}_1}} |\Psi''(\tau)| > T^{5/2}\log T\theta^{-4}\right),$$

$$\mathbf{P}\left(\sup_{\mathcal{V}_{\mathcal{A}_1}} |\mathrm{X}^{(3)}(\tau)| > T^{7/2}\log T\mathcal{M}_T\theta^{-6}\right), \qquad \mathbf{P}\left(\sup_{\mathcal{V}_{\mathcal{A}_1}} |\Psi^{(3)}(\tau)| > T^3\log T\|\lambda^{(3)}\|\theta^{-6}\right)$$

*are negligible in the sense of Definition 1. Moreover,*

$$\mathbf{E}\left(\sup_{\mathcal{V}_{\mathcal{A}_1}} |\mathrm{X}^{(3)}(\tau)|^2\right) = O(T^8\theta^{-12}), \qquad \mathbf{E}\left(\sup_{\mathcal{V}_{\mathcal{A}_1}} |\Psi^{(3)}(\tau)|^2\right) = O(T^8\theta^{-12}\log^{1/2} T).$$

**Proof.** The proof is standard, using Lemmas 14 and 5, and is thus omitted. □

### 4.3. Conclusion of the proof

**Lemma 10.** *Let $\mu_T = T^{1/16}$ and $\mathcal{B} = B(\theta, \theta T^{-1}\mu_T^{-1})$. Then, for any $\tau \in \mathcal{B}$, we have*

$$\Gamma(\tau) - \Gamma(\theta) = -(\tau - \theta)^2 \mathrm{I}(\lambda) + o\{(\tau - \theta)^2 \mathrm{I}(\lambda) + (\tau - \theta) \mathrm{I}(\lambda)^{1/2}\}, \tag{40}$$

$$\mathrm{X}(\tau) - \mathrm{X}(\theta) = 2(\tau - \theta)\{1 + o(1)\} \mathrm{I}(\lambda^2)^{1/2} \mathcal{N} + (\tau - \theta)\mathcal{R}_1(\tau), \tag{41}$$

$$\Psi(\tau) - \Psi(\theta) = (\tau - \theta)\mathcal{R}_2(\tau), \tag{42}$$

*where $\mathcal{N}$ is an $\mathcal{N}(0,1)$ random variable and for $i = 1,2$, the process $\mathcal{R}_i$ satisfies*

$$\exists C_1, C_2 > 0, \forall x \in [0, \log T] \qquad \mathbf{P}\left(\sup_{\mathcal{B}} |\mathcal{R}_i| > x\mathrm{I}_T^{1/2}\right) \leq (1 + C_1\mu_T x)\exp(-C_2 x^2).$$

**Proof.** To prove (40), note that

$$\Gamma(\tau) - \Gamma(\theta) = \sum_k \lambda_k T \left[\left|c_k\hat{\phi}(a_{k,k}) + \sum_{l \neq k} c_l\hat{\phi}(a_{k,l})\right|^2 - \left|c_k + \sum_{l \neq k} c_l\hat{\phi}(b_{k,l})\right|^2\right]$$

$$= \sum_k \lambda_k T |c_k|^2 \{|\hat{\phi}(a_{k,k})|^2 - 1\}$$

$$\quad + \sum_k \lambda_k T \left[\left|\sum_{l \neq k} c_l\hat{\phi}(a_{k,l})\right|^2 - \left|\sum_{l \neq k} c_l\hat{\phi}(b_{k,l})\right|^2\right] + \gamma$$

$$= \alpha + \beta + \gamma,$$

where $\gamma$ regroups the crossed terms coming from each of the squares. But for all real $u$, $\hat{\phi}(u) = 1 - \pi^2 u^2/6 + u^3 \psi(u)$ with $\psi$ bounded on $\mathbb{R}$, so

$$\alpha = -(T^3 \tau^{-2} \theta^{-2}/12) \sum_k \lambda_k (2\pi k)^2 |c_k|^2 (\tau - \theta)^2 \{1 + kT(\tau^{-1} - \theta^{-1}) O(1)\}.$$

Now, note that by (**F1**), the rate of $I(\lambda)$ is a constant times $T^3 \theta^{-4}$. Using the fact that $\sum k^3 |c_k|^2$ is finite and the fact that $\tau$ lies in $\mathcal{B}$, we obtain $\alpha = -I(\lambda)(\tau - \theta)^2 \{1 + o(1)\}$.

To bound $|\beta|$, we use the inequality $||a|^2 - |b|^2| \le |a - b|(|a| + |b|)$ and (24):

$$\left| \sum_{l \ne k} c_l \{\hat{\phi}(a_{k,l}) - \hat{\phi}(b_{k,l})\} \right| \le \sum_l |c_l| \left[ \sup_{[a_{k,l}, b_{k,l}]} |\hat{\phi}'| \right] |a_{k,k}| \le C(2\theta/T) |k| |1 - \theta/\tau| T/\theta,$$

$$\sum_{l \ne k} |c_l| [|\hat{\phi}(a_{k,l})| + |\hat{\phi}(b_{k,l})|] \le C(2\theta/T).$$

Thus, $|\beta|$ is bounded from above by $C \sum_k \lambda_k |k| |\tau - \theta| \le C N_T^2 |\tau - \theta|$, which is $o((\tau - \theta) I(\lambda)^{1/2})$. The same holds for $\gamma$, similarly.

Note that by means of a Taylor expansion, there can be seen to exist random reals $c$ and $d$ such that:

$$X(\tau) - X(\theta) = X'(\theta)(\tau - \theta) + X''(c)(\tau - \theta)^2/2,$$
$$\Psi(\tau) - \Psi(\theta) = \Psi'(d)(\tau - \theta),$$

where the distribution of $X'(\theta)$ is $\mathcal{N}(0, \mathbf{E}(X'(\theta)^2))$ and $\mathbf{E}(X'(\theta)^2)$ is given by Lemma 7.

Let us now control the deviations of $X''(c)$. It suffices to control the supremum of $X''(\tau)$ for $\tau$ in $\mathcal{B}$. Note that for any $\gamma > 0$,

$$\mathbf{E}(\exp\{2\gamma X''(\tau)\}) = \exp\{2\gamma^2 \mathbf{E}(X''(\tau)^2)\}.$$

Using Lemma 5 and the definition of $\mathcal{B}$, we obtain that, for any $y > 0$ and $\gamma > 0$,

$$\mathbf{P}\left( \sup_{\tau \in \mathcal{B}} |X''(\tau)| > y \right) \le 2 \exp(-\gamma y + C T^5 \theta^{-8} \gamma^2) \left( 1 + \gamma C T^{5/2} \int_{\mathcal{B}} \frac{d\tau}{\tau^6} \right)$$
$$\le 2 \exp(-\gamma y + C T^5 \theta^{-8} \mu^2)(1 + \gamma C T^{5/2} N_T \theta^{-4} \mu_T^{-1}).$$

Letting $y = \mu_T T^{5/2} \theta^{-4} x$ and $\gamma = \eta T^{-5/2} \theta^4 x$ with sufficiently small $\eta$, we obtain the desired bound for $\mathcal{R}_1$.

Finally, we control the deviations of $\Psi'$ on $\mathcal{B}$. For any real $\gamma$ such that $\gamma^{-1} > 8 \sup_p \beta_p^{(1)}$, using (37) and the inequality $-\log(1 - u) \le u + u^2$ for $u < 1/2$, we have

$$\mathbf{E}(\exp\{2\gamma \Psi'(\tau)\}) = \exp\left\{ -\sum_p \log(1 - 4\gamma \beta_p^{(1)})/2 \right\} \le \exp\left\{ 2\gamma \sum_p \beta_p^{(1)} + 8\gamma^2 \sum_p \beta_p^{(1)\,2} \right\}.$$

Then, using (38), $\sum_p \beta_p^{(1)} = 0$. Using (39) and Lemma 5, $\sum_p \beta_p^{(1)\,2} \leq T^2 \tau^{-4} \|\lambda^{(1)}\|^2$. In particular, $\sup_p \beta_p^{(1)} \leq CT\theta^{-2}\|\lambda^{(1)}\|$. To conclude, we apply Lemma 14 with $\gamma = \eta \mathrm{I}_T^{-1/2} x$ and sufficiently small $\eta$. □

**Lemma 11.** *As $T$ tends to $+\infty$, for any integer $p$ and sufficiently large $D$,*

$$\mathbf{P}(|\theta^* - \theta|T^{3/2}\theta^{-2} \geq D\log^{1/2} T) = o(T^{-p}). \tag{43}$$

**Proof.** We proceed in two steps. Recall that from Section 3, we have $\theta^* \in B(\theta, \varepsilon_T)$ on $\mathcal{A}_0$, where $\mathbf{P}(\mathcal{A}_0^c)$ is negligible. We first show that $\mathbf{P}(\theta^* \in B(\theta, \theta T^{-1}\mu_T^{-1})^c)$ is negligible. Note that if $V$ is a neighborhood of $\theta$ such that $]k\theta/\tau[=k$ for all $|k| \leq N_T$ and $\tau$ in $V$ (e.g., $V = B(\theta, \varepsilon_T)$), we have

$$\Gamma(\tau) - \Gamma(\theta) = \sum_k \lambda_k T |c_k|^2 (|\hat{\phi}(a_{k,k})|^2 - 1) + O\left(TN_T \frac{\beta_T}{T}\right).$$

Let us define $\mathcal{C} = B(\theta, \varepsilon_T) \cap B(\theta, \frac{\theta}{T\mu_T})^c$. We next prove that $\mathbf{P}(\theta^* \in \mathcal{C})$ is negligible. Note that for $C = \pi^2/2$, $D = \pi^2$ and any $u \in \mathbb{R}$, we have $|\hat{\phi}(u)|^2 - 1 \leq -Cu^2(D+u^2)^{-1}$, thus for any $\tau \in \mathcal{C}$,

$$\Gamma(\tau) - \Gamma(\theta) \leq -T\sum_k \lambda_k |c_k|^2 \frac{Ck^2}{D\mu_T^2 + k^2} + O(N_T\beta_T)$$

$$\leq -T\mu_T^{-2}C \sum_{|k|\leq \mu_T} \lambda_k k^2 |c_k|^2 + O(N_T\beta_T)$$

$$\leq -CT\mu_T^{-2} \qquad \text{since } \mu_T^2 \beta_T N_T = o(T).$$

Recall that $\theta^*$ is defined by (8) and let $\eta = \mathrm{X} + \Psi$. Then

$$\mathbf{P}(\theta^* \in \mathcal{C}) \leq \mathbf{P}(\theta \notin B(e_T, 1/4)) + \mathbf{P}\left(\mathrm{L}(\theta) \leq \sup_{\mathcal{C}} \mathrm{L}(\tau)\right),$$

$$\mathbf{P}(\theta^* \in \mathcal{C}) - \mathbf{P}(\mathcal{A}_0^c) \leq \mathbf{P}\left(\Gamma(\theta) \leq \sup_{\mathcal{C}} \Gamma(\tau) + 2\sup_\Theta |\eta(\tau)|\right) \leq \mathbf{P}\left(2\sup_\Theta |\eta(\tau)| \geq CT\mu_T^{-2}\right).$$

Using Lemma 6, we conclude that $\mathbf{P}(\theta^* \in B(\theta, \frac{\theta}{T\mu_T})^c)$ is negligible.

Letting $\mathcal{A}' = \{\theta^* \in B(\theta, \frac{\theta}{T\mu_T})\}$ and $E_x = \{\frac{\theta^2}{2T\mu_T} > |\tau - \theta| > x\mathrm{I}_T^{-1/2}\}$ for $x > 0$, we have

$$\mathbf{P}(|\theta^* - \theta|\mathrm{I}_T^{1/2} > x) \leq \mathbf{P}\left(\sup_{|\tau-\theta|>x\mathrm{I}_T^{-1/2}} \mathrm{L}(\tau) \geq \mathrm{L}(\theta)\right) \leq \mathbf{P}\left(\sup_{E_x} \mathrm{L}(\tau) - \mathrm{L}(\theta) \geq 0\right) + \mathbf{P}(\mathcal{A}'^c).$$

Now, write $\mathrm{L}(\tau) - \mathrm{L}(\theta) = \Gamma(\tau) - \Gamma(\theta) + \mathrm{X}(\tau) - \mathrm{X}(\theta) + \Psi(\tau) - \Psi(\theta)$. Note that for sufficiently large $T$, by (40), we have $\Gamma(\tau) - \Gamma(\theta) \leq -(\tau - \theta)^2 \mathrm{I}(\lambda)/2$ for any $\tau \in E_x$. Then using expansions (41) and (42), we obtain

$$\mathbf{P}(|\theta^* - \theta|\mathrm{I}_T^{1/2} > x)$$
$$\leq \mathbf{P}\left(\sup_{E_x}(\tau - \theta)[-\tfrac{1}{2}(\tau - \theta)\,\mathrm{I}(\lambda) + 2\,\mathrm{I}(\lambda^2)^{1/2}\mathcal{N} + (\mathcal{R}_1 + \mathcal{R}_2)] \geq 0\right) + \mathbf{P}(\mathcal{A}'^c)$$
$$\leq \mathbf{P}(x \leq C\mathcal{N}) + \mathbf{P}\left(x \leq \mathrm{I}_T^{-1/2}\sup_{E_x}|\mathcal{R}_1 + \mathcal{R}_2|\right) + \mathbf{P}(\mathcal{A}'^c).$$

To conclude, we use the standard bound for the tail of a Gaussian random variable and Lemma 10, and we set $x = D\log^{1/2} T$ with sufficiently large $D$. $\square$

**Lemma 12.** *As $T$ tends to $+\infty$,*

$$\mathbf{E}((\widehat{\tau} - \theta)^2\mathrm{I}_T) = 1 + \{1 + o(1)\}\frac{R_T(f, \lambda)}{\|f'\|^2}.$$

**Proof.** The definition of $\widehat{\tau}$ implies that $\mathbf{E}((\widehat{\tau} - \theta)^2\mathrm{I}_T)[\mathbf{E}(\mathrm{L}''(\theta))]^2 = \mathbf{E}(\mathrm{L}'(\theta)^2)\mathrm{I}_T$. To compute the expectations, note that (38) is zero for $p \geq 1$ and that $\mathrm{X}'(\theta)\Psi'(\theta)$ has zero mean since it is a product of an uneven number of stochastic integrals with respect to $dW(t)$. Thus, $4[\mathbf{E}(\mathrm{L}''(\theta))]^2 = \Gamma''(\theta)^2$ and $4\mathbf{E}(\mathrm{L}'(\theta)^2) = \Gamma'(\theta)^2 + \mathbf{E}(\mathrm{X}'(\theta)^2) + \mathbf{E}(\Psi'(\theta)^2)$.

Using Lemmas 7, 8 and 2, we can write the following expansion of the risk for $\widehat{\tau}$:

$$\mathbf{E}((\widehat{\tau} - \theta)^2) = \left[\mathrm{I}(\lambda^2) + \frac{T^3}{\theta^4}o\left(\frac{\|\lambda'\|^2}{T}\right) + \frac{T^2}{\theta^4}\frac{\pi^2}{3}\sum_k \lambda_k^2 k^2\{1 + o(1)\}\right]$$
$$\times \left[\mathrm{I}(\lambda) + \frac{T^3}{\theta^4}o\left(\frac{\|\lambda'\|^2}{T}\right)\right]^{-2},$$

$$\mathbf{E}((\widehat{\tau} - \theta)^2\mathrm{I}_T) = \left[1 + \|f'\|^{-2}\left\{\sum_k(\lambda_k^2 - 1)(2\pi k)^2|c_k|^2 + \frac{1}{T}\sum_k \lambda_k^2(2\pi k)^2 + o\left(\frac{\|\lambda'\|^2}{T}\right)\right\}\right]$$
$$\times \left[1 + \|f'\|^{-2}\sum_k(\lambda_k - 1)(2\pi k)^2|c_k|^2 + o\left(\frac{\|\lambda'\|^2}{T}\right)\right]^{-2}.$$

Using (11), we can expand the denominator. Then using the expansion of $(1 - \varepsilon)^{-2}$ around $\varepsilon = 0$ and the fact that, due to (**T**), $\sum_k(1 - \lambda_k)(2\pi k)^2|c_k|^2 = o(R_T(f, \lambda)^{1/2})$, we obtain

$$\mathbf{E}((\widehat{\tau} - \theta)^2\mathrm{I}_T) = \left[1 + \|f'\|^{-2}\left\{\sum_k(\lambda_k^2 - 1)(2\pi k)^2|c_k|^2 + \frac{1}{T}\sum_k \lambda_k^2(2\pi k)^2 + o\left(\frac{\|\lambda'\|^2}{T}\right)\right\}\right]$$
$$\times \left[1 - 2\|f'\|^{-2}\sum_k(\lambda_k - 1)(2\pi k)^2|c_k|^2 + o(R_T(f, \lambda))\right],$$

$$\mathbf{E}((\widehat{\tau}-\theta)^2 \mathrm{I}_T) = 1 + \|f'\|^{-2} \sum_k (2\pi k)^2 \left((1-\lambda_k)^2 |c_k|^2 + \frac{1}{T}\lambda_k^2\right) + o(R_T(f,\lambda)).$$

**Lemma 13.** *Let $R_T$ be the functional defined by (12). Then, as $T$ tends to $+\infty$,*

$$\mathbf{E}((\theta^* - \widehat{\tau})^2 \mathrm{I}_T \mathbf{1}_{\mathcal{A}_1}) = o(R_T(f,\lambda)T^{-1}), \tag{44}$$

$$\mathbf{E}((\theta^* - \widehat{\tau})(\widehat{\tau} - \theta)\mathrm{I}_T \mathbf{1}_{\mathcal{A}_1}) = o(R_T(f,\lambda)T^{-1}). \tag{45}$$

**Proof.** Note that by Lemma 1 combined with Lemma 6, the criterion L admits on $\mathcal{A}_0$ a local maximum inside the ball $B(e_T, 1/4)$. Thus, $\mathrm{L}'(\theta^*) = 0$ on this event. Using a Taylor expansion of $\mathrm{L}(\tau)$, there exists $\omega \in [\theta, \theta^*]$ such that

$$\mathrm{L}'(\theta^*) = 0 = \mathrm{L}'(\theta) + (\theta^* - \theta)\mathrm{L}''(\theta) + \tfrac{1}{2}(\theta^* - \theta)^2 \mathrm{L}^{(3)}(\omega).$$

Using (16), we deduce that

$$(\theta^* - \widehat{\tau})\mathbf{E}(\mathrm{L}''(\theta)) = (\theta - \theta^*)\{\mathrm{L}''(\theta) - \mathbf{E}(\mathrm{L}''(\theta))\} - \tfrac{1}{2}(\theta - \theta^*)^2 \mathrm{L}^{(3)}(\omega). \tag{46}$$

Note that using Lemmas 3 and 9,

$$\mathbf{E}\left[\sup_{\mathcal{V}_{\mathcal{A}_1}} |\mathrm{L}^{(3)}(\omega)|^2\right] \leq \mathbf{E}\left[\sup_{\mathcal{V}_{\mathcal{A}_1}} |\Gamma^{(3)}|^2 + |\mathrm{X}^{(3)}|^2 + |\Psi^{(3)}|^2\right] \leq CT^8 \theta^{-12} \log^{1/2} T.$$

Note that on $\mathcal{A}_1$, $(\theta^* - \theta)^2 \leq D^2 \theta^4 T^{-3} \log T$. Moreover, thanks to Lemma 2, $\mathbf{E}(\mathrm{L}''(\theta))^{-2} = 4\mathbf{E}(\Gamma''(\theta))^{-2} \leq (CT^3 \theta^{-4})^{-2}$. Thanks to Lemma 5,

$$\mathbf{E}\{\mathrm{L}''(\theta) - \mathbf{E}(\mathrm{L}''(\theta))\}^2 \leq \mathbf{E}(\mathrm{X}''(\theta)^2) + \mathbf{E}(\Psi''(\theta)^2) \leq CT^5 \theta^{-8}.$$

From the preceding inequalities, we conclude that

$$\mathbf{E}((\theta^* - \widehat{\tau})^2 \mathrm{I}_T \mathbf{1}_{\mathcal{A}_1}) \leq C(\log T)^{5/2} T^{-1} \leq C(\log T)^{-3/2} T^{-1} \|\lambda'\|^2 = o(R_T(f,\lambda)),$$

by (**W1**), which proves (44). Finally, (45) is proved by similar arguments. □

**Lemma 14.** *Let $(\mathrm{X}_t)$ be a stochastic process differentiable a.s., $\mu$ and $x$ positive real numbers and $I$ an interval of $\mathbb{R}$. We then have*

$$\mathbf{P}\left(\sup_{\tau \in I} \mathrm{X}_\tau > x\right) \leq \exp(-\mu x) \sup_{\tau \in I}(\mathbf{E}\exp(2\mu \mathrm{X}_\tau))^{1/2} \left(1 + \mu \int_{\tau \in I}(\mathbf{E}|\mathrm{X}'_\tau|^2)^{1/2} \,\mathrm{d}\tau\right).$$

**Proof.** The proof of this lemma can be found in Golubev [9]. □

## Acknowledgements

The author would like to thank the editor, the associate editor and the referees for their helpful comments.